\documentclass[titlepage,twoside,12pt]{article}
\usepackage{amsmath,amsthm,amssymb,amsfonts,amscd,pstricks}
\begin{document}

{\bf \Large Difficulties with Learning and \\ \\ Teaching Calculus} \\ \\

{\bf Elem\'{e}r E Rosinger} \\
Department of Mathematics \\
and Applied Mathematics \\
University of Pretoria \\
Pretoria \\
0002 South Africa \\
eerosinger@hotmail.com \\ \\

{\bf Abstract} \\

Several thoughts are presented on the long ongoing difficulties both students and academics
face related to Calculus 101. Some of these thoughts may have a more general interest. \\ \\

{\bf 1. What Is Going On ... } \\

For several decades by now around the world, ever new masses of first year students, and not
only in science and engineering, are faced with having to learn what is usually called
"Calculus 101". And in spite - or is it perhaps because of - the ever ongoing "reforms" in
teaching that subject, there seems to be no light at the end of the tunnel ... \\
In fact, one is not so sure whether all of us, students and teachers, are in any sort of
tunnel at all, or rather, we face an insurmountable obstacle ... \\

That long ongoing and widespread situation of deadlock - or is it rather a dead end - has so
far led to one reaction only :

\begin{itemize}

\item Place the whole burden of it on the Departments of Mathematics and the respective
academics, who alone are supposed to be both part of the problem, and of its much sought after
solution.

\end{itemize}

The fact, however, is that one need not be a specialist in Systems Theory or in Control Theory,
and instead, one need only be aware of elementary concepts of management, to realize that the
situation has not been approached in a proper enough manner, even if for some decades by now
so many have been concerned about and involved in it. \\

Here we shall suggest two aspects which are critically important, yet have hardly been given
any attention. Formulated briefly, these aspects are :

\begin{itemize}

\item Seriously incomplete consideration of the situation, with the consequent limitation to
the resulting unilateral actions, actions which cannot but turn out to be unsuccessful.

\item The deeply vulnerable nature of modern technological societies in not given enough
consideration.

\end{itemize}

In this regard it is important to realize that the respective reasons why these two aspects
have been disregarded, or simply missed, are rather natural in our times. Consequently, it is
hard to blame anybody in particular for not giving due attention to them. \\

Therefore, in order to deal with the difficulties in learning and teaching Calculus 101, and
do so at last successfully, we should first go significantly beyond some of the ways of
thinking which happen to prevail at present. \\ \\

{\bf 2. Input-Output System} \\

Learning and teaching Calculus 101 at college or university, or for that matter, learning and
teaching any academic subject, is in fact a classical input-output system :

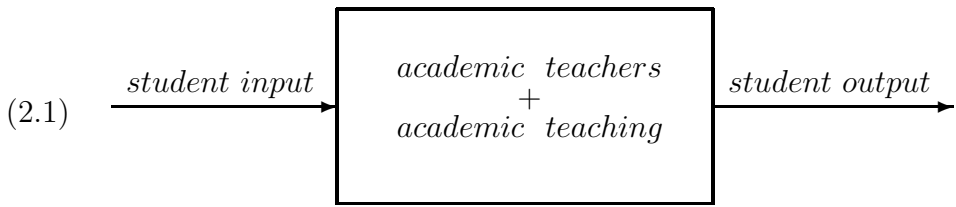
\begin{math}
\setlength{\unitlength}{0.2cm}
\thicklines
\begin{picture}(60,20)

\put(8,10){$student~input$}
\put(7,9){\vector(1,0){15}}
\put(22,15.5){\line(1,0){25}}
\put(0,8){$(2.1)$}
\put(22,15.5){\line(0,-1){13}}
\put(47,15.5){\line(0,-1){13}}
\put(22,2.6){\line(1,0){25}}
\put(26,11){$academic~~teachers$}
\put(34,9){$+$}
\put(26,7){$academic~~teaching$}
\put(47,9){\vector(1,0){16}}
\put(48,10){$student~output$}

\end{picture}
\end{math}

Consequently, there are {\it two factors} which determine the outcome of the process described
by such an input-output system, namely :

\begin{itemize}

\item the quality of the students at "input", and

\item the quality of the academic teachers and academic teaching which is the "transfer
function" in the given input-output system.

\end{itemize}

And there is obviously no way in which to guarantee a satisfactory quality of the students at
"output", without first securing a satisfactory quality of the students at "input". After all,
ever since we have computers, everybody knows the adage : \\

"Garbage in, garbage out." \\

In other words, no matter how good a computer is, if the input data is garbage, the output
data will be quite the same ... \\

The difference with learning and teaching an academic subject at a college or university is
that, indeed, a lot can be done to improve the quality of the academic teachers and academic
teaching. And improving that quality is a rather permanent venture. \\

However, by focusing {\it exclusively} on that issue alone, and not according a comparative
attention to the issue of the quality of the students at "input" can only lead to the
perpetuation of the long ongoing present unsatisfactory situation. \\

Not to mention that, world wide, there does not seem to be any improvement in the quality of
the students at "input", at least not when it comes to Mathematics, and when one considers the
masses of new and new students who are supposed to learn Calculus 101. In other words, the
primary and secondary school system which is supposed to deliver those masses of students at
"input" is simply not able to prepare Mathematically apt students, and do so anywhere near to
the numbers which will have to learn Calculus 101. \\

And as things stand at present, it may be that 10 percent, if not more, of those who complete
the school system may end up having to face Calculus 101 ... \\

Yet in no part of the world, in no nation on the Earth, has ever been made a through enough
study, a study undertaken over a longer time period, about the percentage of those in the
general population who are able to finish successfully a course of Calculus 101. \\

On the other hand, when it comes to professional sport, in more developed countries we have
very good statistics about how many in the general population can become, say, heavy weight
box champions, football, rugby, baseball, basketball, etc., stars, or top athletes in running,
jumping, swimming, and so on. And no matter how lucrative such natural gifts may be to those
who happen to have them, no one is bothered much by the fact that, actually, only a tiny
minority of the general population are able to make it to the respective categories of
professional sport. \\

Perhaps - and hopefully - related to our human ability of learning Calculus 101 we may be
better off than in making it to top levels in sport. \\
However, in our knowledge societies, where technology changes so fast, and it depends so much
on a fast developing science, which on its turn, depends essentially on a Mathematics far more
difficult than mere Calculus 101, it may be high time to establish whether, indeed, 10 percent,
or perhaps more, and who knows, may be less, of the general population could successfully
learn at least Calculus 101 ... \\

The right to reach the top in sport in more developed societies is one of the many rights
accorded to the general population. Just as it is the right to reach the top in education. \\
However, just as in sport, so in education, and in particular, in being able to learn Calculus
101, right does not automatically mean as well the individual's ability to use it and also
benefit from it. \\
Indeed, the right to education, just like the right to sport, is an issue of equal opportunity,
rather than of equal outcome. \\

Fortunately, our modern societies - except for entertainment - depend very little on the
number of top people in sport. \\
On the other hand, the ways we organize our technology nowadays seem to depend on a large
number of young people who should know Calculus 101, not to mention other yet more difficult
Mathematics. \\
And lacking absolutely any systematic study about the capacity of human societies in general
to produce young people able to learn all that Mathematics, we are simply fighting a war in
which we have a seriously deficient idea about our strengths, and our weaknesses ... \\

And then, the "whipping boy" becomes the academic Mathematician and his or her teaching. \\

Yet, it should by now be more than obvious that no matter how much the respective academic
teachers and academic teaching may improve, that may not affect essentially the outcome, when
it comes to the number and quality of students who managed to learn Calculus 101. And the fact
that such an outcome depends at least as much on the quality of the masses of students who
come to be the "input" in that process can only be further disregarded, or even denied, only
at the cost of perpetuating the present unsatisfactory situation ... \\

Several facts and misconceptions may be appropriate to note here. \\

The myth that Mathematics is difficult for most us humans is both true and false. \\
And that it is false has been amply proven throughout history by the fact that nearly every
human, no matter how uneducated otherwise, and no matter whether illiterate, has always been
not only most eager to know the basic Arithmetics involved in counting his or her possessions,
among others, money for instance, but has also succeeded with such a counting. And as we all
know it, and is still clearly shown nowadays in those parts of the world where there are
significant numbers of uneducated and illiterate people, nearly everybody can count money. \\
On the other hand, it is true that Mathematics is difficult for most of the people. Indeed,
each of us can remember how most of one's school colleagues felt a manifest dislike, if not in
fact horror, of Mathematics. \\

And then, the question arises : is Mathematics indeed easy or difficult for us humans ? \\

What is quite clear in this regard is the following. \\

What may make Mathematics difficult for so many is not so much lack of intelligence, as rather
deficiencies in personality. Certainly, as presented during the school years, Mathematics
requires a {\it continuous} dedication and work, since it is built up step by step on all of
its previously taught parts. Thus one cannot so easily be bad at Mathematics, say, in grades 4
to 6, and then suddenly emerge as a star in the higher grades. And unless one is really
interested and likes the subject, it is most unlikely that one would dedicate to it the
sustained effort needed during most of one's childhood in order to avoid having rather fatal
gaps in the subject. \\

Another difficulty with Mathematics is that nowadays there are preciously few school teachers
who are good at it, and who do their best to interest in the subject as many children as
possible. And in the case of such a rather abstract subject like Mathematics, it is crucial
for a teacher to be able to make it liked by as many as possible of one's pupils, since when
left all on its own, Mathematics may easily appear to be a most strange, uninteresting and
pretty irrelevant subject. \\
Consequently, many of those children who would in fact be able to learn Mathematics, are put
off by the subject, and thus end up with gaps which simply cannot ever after be overcome
properly. \\

Related to this latter issue, one should further note that, unlike in earlier times, say,
prior to World War II, there is no longer any particular prestige, let alone pecuniary reward,
in being a school teacher. Added to that comes the fact that those who may be able to learn
enough Mathematics in order to become school teachers can easily find far more convenient and
lucrative jobs, for instance, in the IT industry. \\

As a consequence, the ongoing giant enterprize of teaching Calculus 101 to new and new masses
of students happens to be placed on the proverbial "feet of clay" ... \\
And the unavoidable negative effects of that unfortunate situation are bound to remain with us
for evermore, no matter how much one would flog the hapless academic mathematicians, and no
matter how much they would keep improving their teaching methods ... \\

Let us just remember that in order to get a top class person in sport, one does not start his
or her training when such a person enters college or university. \\
On the other hand, in view of the poor performance related to Mathematics exhibited by the
vast majority of primary and secondary schools, the moment a student faces Calculus 101 is
most likely the first time such a student faces Mathematics in a proper professional
context ... \\

And yet, what should remain as an important and encouraging memento is that, throughout human
history, nearly everybody among no matter how uneducated and illiterate people has been both
eager and able to learn the Mathematics needed for counting one's possessions or
money ... \\ \\

{\bf 3. Increasing the Individual's Insecurity} \\

There are, of course, a number of rather unavoidable reasons why in our modern technological
societies so many of the young people are placed in the situation of having to learn Calculus
101. \\
And some of such reasons may indeed be reasonable, while other ones could possibly be less
so ... \\

Among the latter may simply be those caused by certain over-reactions. \\

One of such cases, long forgotten by so many, happened back in 1958, when following the
successful launching of the first Soviet sputnik in late 1957, the Americans got to feel
deeply shocked and threatened by what they perceived at the time to be the so called "missile
gap". As a consequence, in early 1958, the Eisenhower Administration decided massively to
increase the number of students who would learn science and technology at colleges and
universities. \\

Needless to say, and so unfortunately, the corresponding explosion of student numbers, all of
them having of course to start by studying Calculus 101, was in no way accompanied by a much
needed similarly dramatic improvement of teaching and learning Mathematics across the primary
and secondary schools in America. And in fact, a contrary trend of ever decreasing academic
standards and performance started to prevail in such schools, and not only related to
Mathematics, and not only in America. \\
And this is, among others, how we ended up with the "feet of clay" upon which the teaching of
Calculus 101 to masses of students at colleges and universities has been attempting to stand
for several decades by now ... \\

One effect, not necessarily fortunate, is quite obvious in many places. \\
In earlier times, when in tertiary education only those learned Mathematics who were really
interested in it, learning and teaching that subject was not an issue. And for the vast
majority of students, that is, those not interested in Mathematics, the subject remained
strange, undesired, somewhere far outside of their own world, and definitely in no way
affecting anything at all in their lives ... \\

Nowadays, with so many young people having to face Calculus 101, yet hardly, if at all,
managing to do so, what happens is that we are continually increasing the number of those in
the general population who in their young adulthood had to face a really tough and highly
unpleasant intellectual test, and then failed it, or at best, somehow managed to pass it, but
are fully aware of the fact that they are nowhere near to really master it. \\

Being in such a sort of ever growing category of "de facto intellectually underprivileged" -
or to use the more "politically correct" term, of "de facto intellectually challenged" - is
not quite a joke, or something one can easily forget or disregard. Indeed, it is one thing
never having made it to, say, the quarterback position in the school, college or university
football team, while it is far different - and worse - never having been able to really
understand and master even Calculus 101 ... \\

And to the extent that such a person pursues in his or her life a career in science or
technology, he or she will for ever after remain with a nagging sense of professional
insufficiency, and thus, with a certain amount of intellectual insecurity ... \\

Such a lingering feeling of intellectual insecurity on the part of an ever growing number of
people among the general population may of course have a variety of effects, some of them
perhaps positive, and many other ones negative ... \\

So far, however, no one seems to have given any more systematic consideration to that
issue ... \\
No one seems to have seriously asked the question :

\begin{quote}

"What is the point in having new and new masses of young people pushed into that sort of
lingering feeling of intellectual insecurity ?"

\end{quote}

{~} \\

{\bf 4. Vulnerable Modern Societies} \\

Let us remember that modern technological societies started less than 250 year back, with the
invention of the steam engine in the late 1700s. \\
This watershed event, however, is not given due consideration when it comes to critically
important aspects related to teaching and learning modern science and technology. \\

Indeed, prior to our modern technological era, and for millennia, most of what is called Gross
National Product, or GNP, was produced in agriculture which could easily involve as many as 90
percent of the general population. And all of that agriculture was so primitive, simple and
routine that no one had to give any attention to the next generation of peasants learning what
they had to know. \\

Needless to say, when considered in itself, and certainly not from the point of view of the
masses of peasants involved, that system had several advantages. Among them :

\begin{itemize}

\item It cost society as such next to nothing - both in resources and in training institutions
with qualified personnel - to perpetuate the knowledge it needed for running its production.

\item It was extremely resilient, since in order to destroy such a production system, one had
to destroy the vast majority of the population.

\end{itemize}

Our modern technological societies, on the other hand, quite dramatically {\it lack} both of
these advantages. \\
And that fact has serious consequences which have so far not been appreciated, and thus acted
upon accordingly. It appears that the time passed since the invention of the steam engine was
not long enough for waking up to the radical novelty - and consequent {\it vulnerability} - of
our modern technological societies ... \\

One of the unprecedented and critical aspects of modern technological societies is the
relatively minute number of those who can really master the state of the art aspects of
science and technology, let alone, are able to open up in them new avenues of genuine
importance. \\
And needless to say, the whole of our modern technological societies stand or collapse upon
the existence of such a minute number of people ... \\

Consequently, the timely identification, selection, training and promotion of such people is
of an equally critical importance. \\

In earlier times, it was most likely that one's parents were peasants. And someone whose
parents were peasants was most likely to remain one himself or herself, that sort of rather
automatic and large scale mass process of generational succession not requiring absolutely any
special effort, expense, organization, knowledge, or whatever else of value. And as long as
one had a good enough body physically, one was good enough to become a peasant, since the
knowledge one needed in production could be acquired simply by seeing, doing, and thus
learning without any special organization, during one's childhood. And learn they did in mass
how to do those simple, primitive and routine agricultural tasks of those times. \\

So far, during the last two centuries, it just happened that in what is at present the
developed world the mentioned critically - in fact, vitally - important minute number of
people who can master state of the art science and technology and can further bring it forward
could somehow be obtained. \\
At a closer consideration, however, it is quite obvious that, to a certain not insignificant
extent, and even in the most free societies, the emergence of such people still happens rather
"against the system", than according to any well thought out, well organized, and well
maintained social effort ... \\

Certainly, the way such people manage to emerge cannot be compared with the far more serious
effort, let alone resources, which are invested in producing top level people in sport ... \\

And let us remember what the ancient Romans already knew quite well :

\begin{quote}

"People need bread and circus ..."

\end{quote}

And "bread" comes, of course, before "circus" ... \\

Yet we still care far more about securing those who can deliver the "circus" of sport, than we
do about those who would lead us to more and more "bread" ... \\

And as it happens lately, young people in the Western world go less and less for learning
science and technology. They seem to be, even if intuitively only, well ahead of the rest of
Western societies in realizing the lack of general understanding and appreciation in that part
of the world of the critical, vital in fact, importance of that tiny minority of humans who
can master state of the art science and technology, and can even further open them up to major
new conquest in the future. \\

Regarding the whole of humankind, fortunately, in large societies like India and China, there
is a growing premium on learning science and technology. \\ \\

{\bf 5. Conclusions} \\

Learning Calculus 101 is as much an issue as is teaching it. Focusing alone on improving its
teaching is not going to solve the difficulties we encounter with the massive new and new
waves of students who have to face that subject. \\

It may well happen that, just as with other special human abilities, such as for instance in
art, music or sport, there is in larger normal human societies an upper limit on the number of
those who can successfully deal with Calculus 101, let alone with further yet more difficult
subjects in Mathematics, subjects which nevertheless are essential for state of the art
science and technology. \\
This possibly existing upper limit should be seriously studied. And in case it happens to be
well below the numbers of young people who are presently required to learn Calculus 101, then
corresponding shifts in general policies of modern technological societies should be made. \\
Here one should also note that, as humankind, we only know about Calculus for not much longer
than a mere 300 years. This is certainly rather negligible when compared with many chapters of
Elementary, that is, school Mathematics, which have been known for millennia. It may therefore
be rather natural that, at present, so few of us humans can learn Claculus, let alone learn
the yet more modern and abstract branches of Mathematics. \\
Such a limitation, however, need not necessarily be put on the account of the human species as
a species, since it may actually be a temporary one, even if not one that may diminish
significantly in just about a few more generations. \\

The right to learning, which is a welcome modern development in many societies, cannot be
identified with one's individual capability to do so, no matter how much desirable such an
identification would appear to be. \\

In view of the vital role the tiny minority of humans who can master the state of the art
science and technology, and moreover, can open up horizons for further major discoveries,
should be the object of a general attention and care not less than is the case with the
identification, selection, training and promotion of top professionals in sport. \\

And such a care cannot simply be limited to "throwing more money" at science and technology.
Instead, it requires a system which so far in human history has never had a precedent in its
sophistication. \\
What is done at present, and has historically been done in this regard by universities and
research institutes only corresponds to the ad-hoc, haphazard, artisan sort of approach of the
issues involved. And unfortunately, we are not yet over the deeper and longer lasting negative
effects of the ways the sudden explosion in the number of researchers in science and
technology happened in the 1960s got managed. In those times, due to the ongoing "Cold War",
no one seemed to have the respite to think more deeply and fully about all that would be
involved in such an explosion. Ever since, the system established in the 1960s has been left
to function, whether it was right or wrong, whether it went the right way, or the wrong one.
And all its many problems have been, and are still considered to be "only a problem of more
money" ... \\
Nowadays, when we are embarked upon a longer lasting "War on Terror", it may appear equally
unlikely that our societies may be able to do anything else but let the old, 1960s system run
according to its own inadequate logic, and on rare occasions, perhaps "throw some more money"
at it ... \\

The overall feeling of those among us who have for longer been invovled in scientific research
in one or another of the fields of what is usually called "hard science", is that :

\begin{quote}

"Science is not done scientifically."

\end{quote}

And as things stand nowadays, the aspiration to ever have science done more scientifically
seems not to have its proper time for fulfillment ...

\end{document}